\documentclass{iopart}
\usepackage[ngerman,english]{babel}
\usepackage{iopams,amssymb,amsthm}
\usepackage{bm}
\usepackage{graphicx}

\eqnobysec
\newtheorem{thm}{Theorem}

\newtheorem{rem}[thm]{Remark}

\newcommand{\R}{{\mathbb R}}

\newcommand{\argmin}{\mathop{\rm argmin}}%

\newcommand{\bE}{{\mathbf E}}

\newcommand{\eps}{\epsilon}

\newcommand{\wrt}{with respect to }

\begin{document}

\title{Efficient Rare Event Simulation by Optimal Nonequilibrium Forcing }

\author{Carsten Hartmann, Christof Sch\"utte}
\eads{\mailto{chartman@mi.fu-berlin.de}, \mailto{schuette@mi.fu-berlin.de} }

\date{\today}

\begin{abstract}
Rare event simulation and estimation for systems in equilibrium are among the most challenging topics in molecular dynamics. 
As was shown by Jarzynski and others, nonequilibrium forcing can theoretically be used to obtain equilibrium rare event statistics. The advantage seems to be that the external force can speed up the sampling of the rare events by biasing the equilibrium distribution towards a distribution under which the rare events is no longer rare. Yet algorithmic methods based on Jarzynski's and related results often fail to be efficient because they are based on sampling in path space. We present a new method that replaces the path sampling problem by minimization of a cross-entropy-like functional which boils down to finding the optimal nonequilibrium forcing. We show how to solve the related optimization problem in an efficient way by using an iterative strategy based on milestoning.
\end{abstract}



\section{Introduction}

Molecular dynamics (MD) simulations allow for analysis and understanding of the dynamical behaviour of molecular systems. 
However realistic simulations on timescales beyond microseconds are still infeasible even on the most powerful general purpose computers,
which renders the MD-based analysis of many biological equilibrium processes, that are often rare compared to the characteristic time scale of the system and hence require prohibitively long simulations, impossible. The hallmark of these rare events is that the average waiting time between the events is orders of magnitude longer than the timescale of the switching event itself. Thus rare event simulation and estimation are among the most challenging topics in molecular dynamics. 

The molecular dynamics literature on rare event simulations is rich. Since direct numerical equilibrium simulation is infeasible, all available techniques try to sample from the rare event statistics by biasing the system in one or the other way. Roughly speaking, we can distinguish between two major classes of sampling techniques: class $A$ consists of splitting methods that decompose state space, but are still essentially based on an equilibrium distribution, whereas methods from class $B$ proceed by driving the system under consideration into a nonequilibrium regime that changes the rare events statistics. For a general overview of Monte-Carlo methods for rare events in other application fields, we refer to the textbook \cite{bucklew2004}. 

The list of methods in class $A$ range from reaction-coordinate based techniques via path-space oriented techniques to approaches based on interface sampling or generalized dynamics.  Reaction-coordinate based techniques consider the marginal of the equilibrium distribution in some low-dimensional collective variables  like in direct free energy calculations \cite{freeE}; they suffer from the fact that appropriate reaction coordinates are often not available. Path-space oriented techniques approximate the most important reaction paths that govern the rare event statistics either by sampling distribution of reactive paths like in transition path sampling (TPS) \cite{TPS,chandler1998} or by optimizing an appropriate path functional like in the string method \cite{string}; they become problematic if the path space distribution is multi-modal or generally too complex (e.g., involving bifurcations). Interface sampling techniques like milestoning \cite{Elber} or forward flux sampling (FFS) \cite{FFS} place a set of suitably chosen interfaces in state space between the initial and final state and use them to follow the transition of the system in an iterative manner using equilibrium trajectories that connect neighbouring interfaces. The idea of \emph{generalized dynamics} such as hyperdynamics \cite{voter1997}, metadynamics \cite{laio2002}, conformational flooding \cite{grubmueller1995}, or the adaptive biasing force (ABF) method \cite{ABF} is to bias the system on-the-fly (e.g., by filling in certain energy wells in which the system got trapped during a simulation) so as to enhance rare transitions between metastable states. Although seemingly different, generalized dynamics belong to class $A$, in that they only alter the underlying equilibrium distribution along a predefined set of low-dimensional collective variables. Although these methods have proven to be very efficient, they require that the interesting processes can be described by a few collective coordinates that have to be known in advance. 

Class $B$ consists of methods based on the Jarzynski and Crooks formulae \cite{jarzynski1997,crooks1998} that relate the equilibrium Helmholtz free energy to the nonequilibrium work exerted under external forcing. Instances of nonequilibrium simulations that mimic experiments on controlling and manipulating single molecules (see, e.g., \cite{rief1997,merkel1999}) are single-molecule pulling \cite{hummer2001}, steered molecular dynamics \cite{schulten2004} or bridge sampling \cite{minh2009}, to mention just a few. The corresponding path functionals have the form of cumulant-generating functions for the exerted work \cite{kurchan1998,lebowitz1999} which poses immense challenges to Monte-Carlo simulations and limits the usability of the formulae in practice. Roughly speaking, the usability is limited by the fact that the likelihood ratio between equilibrium and nonequilibrium trajectories is highly degenerate, for the overwhelming majority of nonequilibrium forcings generate trajectories that have almost zero weight \wrt the equilibrium distribution that is relevant for the rare event; cf.~also the discussion in \cite{lelievre2010}. Nevertheless the underlying idea is appealing and a cleverly designed external force may speed up the sampling of the rare events by biasing the equilibrium distribution of the system towards a distribution under which the rare events is no longer rare, while giving numerical estimators that are useful in terms of variance and convergence properties.    

The method presented in this article belongs to the latter class, but shares somes ideas with ideas from class $A$. It takes up the idea that external forcings can speed up the rare event but avoids sampling issues related to nonequilibrium processes. Instead it uses optimal nonequilibrium forcing in connection with splitting methods such as FFS or milestoning, in the sense that the new method uses interfaces to follow the transition of an \emph{optimally driven} system where the external forcing that drives the system from one interface to the next results in a considerable speed-up compared to FFS or milestoning. Specifically, the new method replaces the path sampling problem using an exponential change of measure that can be explicitly computed by minimizing a cross-entropy-like functional, which then yields the optimal forcing. Although the minimization involves solving an optimal control problem, the numerical effort can be drastically reduced when the minimization is done in a clever way; one reason is that the path functional becomes linear after the change of measure whereas it was exponential in the original cumulant-generating function. 

Transformations based on exponential change of measures have a rich tradition in the (risk-sensitive) optimal control literature \cite{james1992,daipra1996,fleming2006} and the theory of large deviations \cite{fleming1977,whittle1994}, and are regularly rediscovered---mostly aiming at turning certain optimal control problems into linearly solvable sampling problems \cite{kappen2005,todorov2009,todorov2012}; cf.~also \cite{schuette2012,toussaint2012}. Here we pursue the reversed strategy and turn a difficult rare event estimation problem into an optimal control problem that can be solved by minimizing a suitable functional.  Thus the basic outline of the new method is: iteratively determine the optimal nonequilibrium forcing by an optimization procedure based on milestoning ideas that avoid path-space sampling and compute the equilibrium rare event statistics from the optimal nonequilibrium forcing.  

Besides introducing the new method the purpose of this article is to explain the basic ideas of how to use optimal control for the estimation and simulation of rare events. Therefore we present only the simplest possible scenario (a particle following an overdamped Langevin dynamics in a conservative force field), without paying too much attention to complete generality or mathematical rigour. The first issue in Section \ref{sec:free} then is to introduce the variational characterization of (generalized) free enregy and the exponential change of measures that are the basis of our optimal control approach. The precise formulation of the optimal control problem, a stochastic control problem with quadratic control costs and an indefinite time horizon, is given in Section \ref{sec:ocp}. In Section \ref{sec:milestoning} we describe the numerical method for computing the optimal control, based on an inexact gradient descent in connection with a milestoning algorithm, and apply it to the controlled first passage between metastable sets. We briefly summarize the results in Section \ref{sec:conclusion} and sketch possible generalization that have been omitted for the sake of brevity.

\section{A variational characterization of free energy}\label{sec:free}

We consider a particle with position $X_{t}\in \R^{n}$ at time $t>0$ which moves in an energy landscape $V\colon\R^{n}\to\R$ according to the equation  
\begin{equation}\label{sde}
dX_{t} = -\nabla V(X_{t})dt + \sqrt{2\eps}\,dB_{t}\,, \quad X_{0}=x \,.
\end{equation}
Here $B_{t}$ denotes standard $n$-dimensional Brownian motion, and $\eps>0$ is the temperature of the system. Under mild conditions on the  energy landscape function $V$ we have ergodicity, and the law of $X_{t}$ converges to a unique equilibrium distribution with density 
\begin{equation*}
\rho(x) = Z^{-1} \exp(-\eps^{-1} V(x))\,,\quad Z=\int_{\R^{n}}\exp(-\eps^{-1} V(x))\,dx\,.
\end{equation*}  
We assume throughout that the temperature is small, relative to the largest energy barriers, i.e., $\eps\ll \Delta V_{\rm max}$. As a consequence, the relaxation of the dynamics towards equilibrium is dominated by the rare transitions over the largest energy barriers. 

Let $W$ be a random variable that depends on the sample paths $(X_{t})_{0\le t\le \tau}$ up to a stopping time $\tau$. We will call $W$ \emph{work} in the following. Given some continuous function $f\colon\R^{n}\times[0,\infty)\to\R$, we suppose that it can be expressed as\footnote{The following considerations below are not at all limited to systems of the form (\ref{sde}) and path functionals like (\ref{work}) and can be easily can be easily generalized to, e.g., non-gradient systems with multiplicative and/or degenerate noise or observables $f$ that are explicitly time-dependent.}
\begin{equation}\label{work}
W  = \int_{0}^{\tau}f(X_{t})\,dt\,.
\end{equation}

Let us further denote by $P$ the probability measure on the space of continuous trajectories that is generated by the Brownian motion in (\ref{sde}), and let $\bE^{x}[\cdot]=\bE[\cdot|X_{0}=x]$ be the expectation \wrt $P$, i.e., the average over all realizations of $X_{t}$ starting at $X_{0}=x$. We call the quantity
\begin{equation}\label{cgf}
F(x) = -\eps \log \bE^{x}[\exp(-W/\eps)]\,.
\end{equation}
the \emph{(conditional) free energy} of $W$ \wrt $P$.

\begin{rem}
Clearly, the functions and the expectation on the right hand side of (\ref{cgf}) do not commute, and it follows by Jensen's inequality that $F(x)\le \bE^{x}[W]$, in accordance with the second law of thermodynamics. But $F$ encodes information about the cumulants of the work $W$ (assuming they exist), namely, 
\begin{equation*}
F(x) = \bE^{x}[W] + \frac{1}{2\eps}\bE^{x}\big[(W-\bE^{x}[W])^{2}\big] + \ldots\,.
\end{equation*}
\end{rem}

\begin{rem}
The similarity between (\ref{cgf}) and Jarzynski's formula \cite{jarzynski1997} is no coincidence. If $\tau=T$ is a deterministic stopping time and $W$ is the nonequilibrium work done on a system during a transition between two equilibrium states $E_{1}$ and $E_{2}$, then $F(E_{1})$ equals the equilibrium free energy difference between $E_{1}$ and $E_{2}$.
\end{rem}

The phrases "work" for the quantity $W$ defined in (\ref{work}) and "free energy" for $F$ as of (\ref{cgf}) are just used to relate to Jarzynski's formula. The framework is much more general as the following example will show.

\paragraph{Guiding example. }

One example, of which we will consider variants below, is the first hitting time of a subset of state space. To this end let $S\subset \R^{n}$ a set and define
\[
\tau=\inf\{t>0\colon X_{t}\in S\}\,.
\]
to be the first time at which $X_{t}$ hits $S$. Choosing the constant function $f=\sigma$ in (\ref{work}), the free energy 
\begin{equation*}
F_{\sigma}(x) = -\eps \log \bE^{x}[\exp(-\sigma\tau/\eps)]\,.
\end{equation*}
considered as a function of $\sigma$ is the scaled cumulant-generating function of $\tau$ when $X_{t}$ is started at $X_{0}=x$. In particular, we can compute the mean first hitting time by
\begin{equation*}
\eps\left.\frac{d F_{\sigma}}{d\sigma}\right|_{\sigma=0} = \bE^{x}[\tau]\,.
\end{equation*}

\subsection{Relative entropy and change of measures}

The strict convexity of the exponential function implies that equality $F(x)=\bE^{x}[W]$ is only attained if $W$ is $P$-almost surely constant; one such case is the adiabatic limit  
\begin{equation*}
W  = \lim_{T\to\infty} \frac{1}{T}\int_{0}^{T}f(X_{t})\,dt\,.
\end{equation*} 

We will restore (\ref{cgf}) to an expression that becomes linear in $W$ after a suitable change of measure. To this end let $Q$ denote a probability measure on the space of continuous trajectories that is absolutely continuous \wrt $P$ (i.e., $\varphi=dQ/dP$ exists). We define the relative entropy of $Q$ \wrt $P$ as
\begin{equation}\label{KLdiv}
I(Q\|P) = \int_{\R^{n}}\log\left(\frac{dQ}{dP} \right)dQ\,.
\end{equation}
(This is also called the \emph{Kullback-Leibler divergence}.) 
We declare that $I(Q\|P)=\infty$ if $Q$ is not absolutely continuous \wrt $P$. Then, by Jensen's inequality, 
\begin{equation}\label{equivcgf}
\eqalign{
F(x) & = -\eps \log \bE^{x}[\exp(-W/\eps)]\\
 &= -\eps \log \bE^{x}_{Q}[\exp(-W/\eps - \log\varphi)]\\
& \le \bE^{x}_{Q}[W] + \eps I(Q\|P)\,,}
\end{equation}
where we have used the notation $\bE_{Q}[\cdot]$ to denote the expectation \wrt $Q$. The last inequality that appears in the literature in various forms as \emph{second-law-like identity} or \emph{generalized Jarzysnki inequality} (cf.~\cite{ueda2010,horowitz2010}) suggests that the free energy and the relative entropy are related by a Legendre-type transformation, viz., 
\begin{equation*}
F(x) = \inf_{Q}\left\{\bE^{x}_{Q}[W] + \eps I(Q\|P)\right\}, 
\end{equation*}
and a result in \cite{daipra1996} implies that the infimum exists and is attained when $Q$ runs over all path measures that are absolutely continuous \wrt $P$. 
By the strict convexity of the exponential function, the latter implies that $W+\eps\log\varphi$ is $Q$-almost surely constant. 

The idea of the approach sketched below then is to represent $Q$ in terms of suitable (parametric) control variables and minimize the right hand side of (\ref{equivcgf}) over all admissible controls.

\section{An optimal control problem}\label{sec:ocp}

The aim of this section is to derive necessary and sufficient conditions for the optimal change of measure that turns (\ref{equivcgf}) into an equality. To this end we follow ideas by Fleming and co-workers \cite{fleming1977,dupuis1997} and consider the exponential cost functional:  
\begin{equation}\label{mgf}
\psi(x) = \bE^{x}\left[\exp\bigg(-\eps^{-1} \int_{0}^{\tau}f(X_{s})\,ds\bigg)\right].
\end{equation}
For a stopping time $\tau$ that is the first hitting time of a set $S\subset\R^{n}$, the Feynman-Kac formula \cite{oksendal2003} implies that $\psi$ solves the elliptic boundary value problem
\begin{equation}\label{FK}
\eps L\psi = f\psi \,,\quad \psi|_{\partial S}  = 1\,,
\end{equation}
where 
\begin{equation}\label{L}
L= \eps\,\nabla^{2} + \nabla V\cdot\nabla \,.
\end{equation}  
is the infinitesimal generator of $X_{t}$, defined on a suitable subspace of $L^{2}(\R^{n})$. We want to transform the boundary value problem (\ref{FK}) into an equation for the unknown control variable in (\ref{equivcgf}). For this we proceed in two steps.

\paragraph{Step 1:} We can safely assume that $\tau$ is almost surely finite. As a consequence, the function $\psi$ in (\ref{mgf}) admits a formal representation of the form
\begin{equation*}
\psi = \exp(-F/\eps) \,.
\end{equation*}
We seek an equation for the free-energy $F$. By chain rule, it follows that     
\begin{equation*}
\eps\exp(F/\eps)L \exp(-F/\eps)= -L F + |\nabla F|^{2}\,,
\end{equation*}
which entails that (\ref{FK}) is equivalent to  
\begin{equation}\label{HJB}
L F- |\nabla F|^{2} + f = 0\,,\quad F|_{\partial S} = 0\,,
\end{equation}
The last equation is known as the Hamilton-Jacobi-Bellmann (HJB) equation of optimal control \cite{fleming2006}; its solution is called \emph{value function} or \emph{optimal cost-to-go}. 

\paragraph{Step 2:} To reveal the stochastic optimal control problem that corresponds to the HJB equation (\ref{HJB}), we first note that 
\begin{equation*}
-|\nabla F|^{2} = \min_{c\in\R^{n}}\left\{\sqrt{2} c\cdot\nabla F + \frac{1}{2}|c|^{2}\right\}\,,
\end{equation*}
from which we recognize that (\ref{HJB}) is equivalent to   
\begin{equation}\label{HJBmin}
\min_{c\in\R^{n}}\left\{ L(c) F + g(x,c) \right\} = 0 \,,\quad F|_{\partial S}  = 0\,,
\end{equation}
with the shorthands
\begin{equation*}
g(x,c) = f(x) + \frac{1}{2}|c|^{2}\,
\end{equation*}
and 
\begin{equation*}
L(c) = \eps\,\nabla^{2} + (\sqrt{2} c - \nabla V)\cdot\nabla \,.
\end{equation*}  
Equation (\ref{HJBmin}) is the Hamilton-Jacobi-Bellman equation of the following optimal control problem that should be compared to the right hand side of (\ref{equivcgf}): minimize 
\begin{equation}\label{indefinite1}
I(u) = \bE\left[\int_{0}^{\tau} g(X_{t},u_{t})\,dt\right]
\end{equation}
over an admissible set $U$ of control laws $u$ with values in $\R^{n}$ and subject to the tilted dynamics
\begin{equation}\label{indefinite2}
dX_{t} = \left(\sqrt{2} u_{t} - \nabla V(X_{t})\right)dt + \sqrt{2\eps}\, dB_{t}\,.
\end{equation}
That is, the expectation in (\ref{indefinite1}) has to be taken wrt the path measure $Q$ generated by the dynamics given by (\ref{indefinite2}).  

\begin{rem}
The dynamics that generates the new path measure $Q$ is again of gradient form if $u=u^{*}$ is the optimal Markovian feedback control, i.e. when $Q=Q(u^{*})$. As a consequence, the optimally controlled process satisfies detailed balance \cite{lebowitz1999}. Indeed, since (\ref{indefinite1}) is quadratic and (\ref{indefinite2}) is affine in the control, the minimizer
\begin{equation*}
c^{*}(x) = \argmin_{c}\left\{ L(c)F + g(x,c) \right\}\,,
\end{equation*}
in (\ref{HJBmin}) is unique (provided that $F$ is sufficiently smooth). The optimal feedback law is then given by $u^{*}_{t} = -\sqrt{2}\nabla F(X_{t})$ and gives rise to the tilted dynamics  
\begin{equation*}
dX_{t} = - \nabla G(X_{t})dt +\sqrt{2\eps}\,  dB_{t}\,, \quad X_{t}\in \R^{n}\setminus S\,,
\end{equation*}
with the tilted potential 
\begin{equation*}
G(x) = V(x) + 2F(x)\,.
\end{equation*}

\end{rem}

\paragraph{Guiding example, cont'd.}
In some cases it is helpful to pursue a reverse strategy and transform the nonlinear HJB equations of an optimal control problem into a linear equation that may be easier to solve (cf. \cite{kappen2005,todorov2009}). 

Consider a Brownian particle under a microscope with a moveable object holder. Let $D\subset\R^{2}$ denote the microscope's focal disc, $X_{t}\in\R^{2}$ the particle position at time $t>0$, relative to the position of the object holder, and $u_{t}$ the motor force. The control task is to move the object holder such that the particle stays in the focus as long as possible. Hence the control objective is the maximization of the mean first exit time from $D$ which amounts to minimizing the cost functional 
\begin{equation*}
I(u) = \bE\left[- \tau + \frac{1}{2}\int_{0}^{\tau} |u_{t}|^{2}\,dt\right],
\end{equation*}
subject to
\begin{equation*}
dX_{t} = \sqrt{2} u_{t} + \sqrt{2\eps}\,dB_{t}\,.
\end{equation*}
Let  
\begin{equation*}
F(x) = \min_{u\in U}\bE^{x}\left[- \tau + \frac{1}{2}\int_{0}^{\tau} |u_{t}|^{2}\,dt\right],
\end{equation*}
be the value function (free energy) of the problem and
\begin{equation*}
\psi(x) = \bE^{x}[\exp(\tau/\eps)] \,.
\end{equation*}
Then the linear boundary value problem for $\psi=\exp(-F/\eps)$ is a Helmholtz equation with Dirichlet boundary conditions, 
\begin{equation*}
\eps^{2} \nabla^{2}\psi + \psi =  0 \,,\quad \psi|_{\partial D} = 1 \,,
\end{equation*}
which can be solved by standard means.

\section{Greedy milestoning algorithm}\label{sec:milestoning}

At first sight it seems that we have not gained much, for we have transformed the original path sampling problem into a complicated nonlinear optimal control problem. However the optimal control formulation opens up other options for the numerical treatment of the rare event sampling in terms of a minimization problem. Another advantage is that it is relatively easy to construct unbiased estimators of the control functional, avoiding both bias and variance issues when estimating exponential observables such as (\ref{cgf}).

\paragraph{Discretization}

Together with the information that the optimal Markov control is of feedback form our minimization problem (\ref{indefinite1})--(\ref{indefinite2})
takes the form
\begin{eqnarray*}
F(x) & = & \min_{u_t=c(X_t)} \bE^x_Q\left[\int_{0}^{\tau} g(X_{t},u_{t})\,dt\right]\label{final}
\end{eqnarray*}
with $Q$ denoting the path measure generated by the dynamics given by (\ref{indefinite2}).  
We discretize this optimization problem by choosing a finite dimensional ansatz space for the space of admissable feedback functions $c$:
We choose sufficiently smooth and integrable vector fields $b_j\colon\R^{n}\to\R^{n}$, $j=1,\ldots,m$, so that  
\[
c(x)=\sum_{j=1}^m a_j b_j(x)\,,\quad a_{j}\in \R\,,
\]
or, respectively, we choose scalar ansatz function $v_j\colon\R^{n}\to\R$, $j=1,\ldots,m$, so that
\[
F(x)=\sum_{j=1}^m a_j v_j(x),\quad b_j=-\sqrt{2}\nabla v_j\,.
\]
The minimization problem then amounts to minimizing the cost functional
 \begin{equation}\label{semicost}
\tilde{I}(a) = \bE_{Q}\left[\int_{0}^{\tau}\bigg(f(X_{s}) + \frac{1}{2}\Big|\sum_{j}a_{j}(s)b_{j}(X_{s})\Big|^{2} \bigg)ds\right]
\end{equation}
over the unknown coefficients $a=(a_{1},\ldots,a_{m})$ where $Q=Q(a)$, the path measure of the controlled diffusion (\ref{indefinite2}) also depends on the coefficients; for the moment we remain with the imprecise statement that the measure $Q$ has a density $\varphi(\cdot;a)$ \wrt a (fictitious) uniform measure on the space of all continuous paths in $\R^{n}$, which is a function of the unknown coefficients.\footnote{More precisely, $Q=Q^{\delta}_{x}$ is the probability to find paths $(X_{s})_{0\le s\le T}$ in a small tube around a smooth curve $\gamma:[0,T]\to\R^{n}$, i.e., $Q^{\delta}_{x}(\gamma)=P(\|X_{s}-\gamma(s)\|\le \delta\,|\,X_{0}=x)$. By the Girsanov theorem, $Q_{x}=\lim_{\delta\to 0}Q^{\delta}_{x}$ has a density $\varphi=\exp(-S(\gamma))$ \wrt the Gaussian measure induced by the Brownian motion $\tilde{B}_{s}=x + \sqrt{2\eps}B_{s}$, where $S(\gamma)$ is the Onsager-Machlup functional \cite{duerr1978}.}

\paragraph{Gradient descent} 

We minimize the cost functional $\tilde{I}(a)$ by a doing a gradient descent in the coefficient vector $a=(a_{1},\ldots,a_{m})$. Specifically, we iterate the map 
\begin{equation*}
a^{(i+1)} = a^{(i)} - \alpha_{i} \nabla \tilde{I}\big(a^{(i)}\big),
\end{equation*}
where $i$ is the iteration index and $(\alpha_{i})_{i\ge 1}$ is a bounded sequence of stepsizes for the gradient search. For instance, we can do a line search in the descent direction and determine $\alpha_{i}$ so that it satisfies the Wolfe condition \cite{wright1999}. Details of the iteration that is based on an Euler-Maruyama discretization of the path measure $Q$ will be given below in the appendix. The overall algorithm thus has the following steps:
\begin{itemize}
\item Choose scalar-valued ansatz functions $v_j$ with support in the interesting region of state space and related vector fields $b_j=-\sqrt{2}\nabla v_j$.
\item Choose initial coefficients $a^{(0)}=(a_j^{(0)})$ such that the free energy or value function $\sum_{j=1}^m a_j v_j(x)$ fills up the main wells in the energy landscape $V$.
\item Iterate the following steps in $i$, starting with $i=0$, until a prescribed termination criterium is satisfied: 
\begin{enumerate}
\item Sample the path measure $Q=Q(a^{(i)})$ and evaluate $\nabla \tilde{I}\big(a^{(i)}\big)$ (see formula (\ref{gradI}) in the appendix).
\item Perform a gradient descent $a^{(i+1)} = a^{(i)} - \alpha_{i} \nabla \tilde{I}\big(a^{(i)}\big)$.
\end{enumerate}
\end{itemize}

\begin{rem}
The gradient search algorithm can be regarded as a variant of the cross-entropy method that is a relatively new Monte-Carlo technique for the sampling of rare events which goes back to Rubinstein and others \cite{rubinstein2004}. It is based on the idea that an optimal change of measure can be found by minimizing the Kullback-Leibler divergence (\ref{KLdiv}) over a family of probability measures $Q$ in terms of the tilting parameter $c$. Compared to equilibrium rare event simulation algorithms used in molecular dynamics using the optimal change of measure has the advantage that the likelihood ratio $dQ/dP$ stays of order one, while rare events under the original dynamics (here: diffusion in an energy landscape $V$) are no longer rare under the forced dynamics (\ref{indefinite2}). As a consequence, sampling the path measure $Q$ is significantly more efficient than sampling the original path measure $P$ since the trajectories to be sampled from $Q$ are much shorter on average (i.e., the expected hitting time is considerably shorter).
 \end{rem}

\paragraph{Milestoning algorithm} 

\begin{figure}
 \centering
 \includegraphics[scale = 0.5]{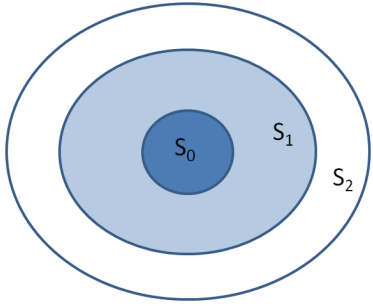}     
 \caption{Illustration of nesting of sets for the milestoning iteration. }
 \label{f:illu}
\end{figure}

For problems with a large state space or for strongly metastable systems, the above algorithm may still be inefficient since sampling the path measure $Q$ may involve many rather long trajectories. In this case the computation can be broken down to transitions between neighbouring interfaces as in milestoning \cite{Elber} or in FFS \cite{FFS}. We explain the basic steps of this procedure: Let 
 \[
\tilde{F}(x) = \min_{a} \bE^{x}_{Q}\left[\int_{0}^{\tau} \tilde{g}(X_{s},c(X_{s})) \,ds\right]
\]
denote the semi-discretized value function of the problem, with the shorthand
 \[
\tilde{g}(x,c(x)) = \sigma f(x) + \frac{1}{2}\Big|\sum_{j}a_{j}b_{j}(x)\Big|^{2}\,.
\]
Suppose that $S=S_0$ is the set of interest and $\tau=\tau_{0}$ is the first hitting time of $S_{0}$; we now choose nested sets or \emph{milestones} $S_{0}\subset S_{1}\subset S_{2}\subset \ldots$ (cf.~Figure \ref{f:illu}). We first compute $\tilde{F}$ in $S_1\setminus S_0$ by finding the optimal control policy $c$ in $S_1\setminus S_0$. That is, our ansatz functions in the above gradient descent algorithm only have to be non-vanishing in $S_1\setminus S_0$. In particular this gives $\tilde{F}$ on $\partial S_1$, the outer boundary of $S_1\setminus S_0$. We can repeat the same algorithm in the set $S_2\setminus S_1$; then letting $x\in  S_2\setminus S_1$ and letting $\tau_1$ denote the first entry time into $S_1$, we have
\[
\tilde{F}(x) =\min_a \bE^{x}_{Q} \left[\int_0^{\tau_1} \tilde{g}(X_s,c(X_s))ds + \tilde{F}(X_{\tau_1}) \right],
\] 
where $X_{\tau_1}\in \partial S_1\subset S_1\setminus S_0$ for which $\tilde{F}$ has been computed in the previous step. By iterating the algorithm we eventually obtain $\tilde{F}$ on all set boundaries $\partial S_{i}$, $i=0,1,2,\ldots$. Thus, the milestoning iteration can be implemented as an outer loop which contains the above gradient descent algorithm in every of its iterations.

 \begin{rem}
 The milestoning variant of the gradient descent algorithm only requires the computation of an ensemble of \emph{short} trajectories of the controlled system (\ref{indefinite2}). Here "short" means that they are orders of magnitude shorter than those in typical path-space sampling algorithms like TPS, and equilibrium milestoning or FFS.
 \end{rem}

\subsection{Guiding example: computing the mean first passage time}

\begin{figure}
 \centering
 \includegraphics[width = 0.65\textwidth]{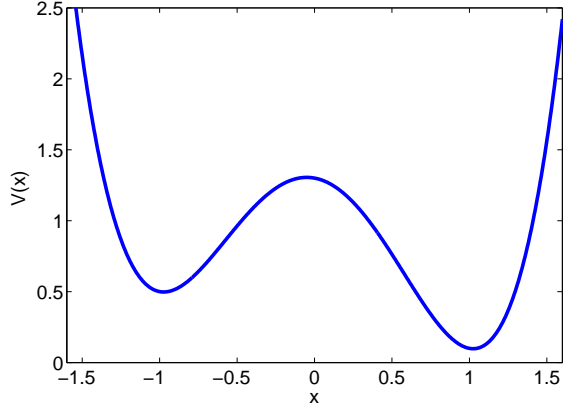}     
 \caption{Skew double-well potential $V$. }
 \label{fig:pot}
\end{figure}

We consider the uncontrolled dynamics (\ref{sde}) with the one-dimensional potential shown in Figure \ref{fig:pot}. Suppose we are interested in computing the  mean first passage time to the set $S=[-1.1,-1]$ in terms of the free energy (\ref{cgf}). Let  
\begin{equation*}
\tau = \inf\{t>0\colon X_{t}\in \partial S\}\,.
\end{equation*} 
be the first hitting time of $S$, consider the constant function $f=\sigma$, and the scaled moment generating function
\[
\psi_{\sigma}(x)=\bE^{x}[\exp\left(-\sigma\tau/\eps\right)],
\]
considered as a function of $\sigma$. The quantity of interest is the mean first passage time of the uncontrolled dynamics, 
\[
-\eps\left.\frac{d\psi_{\sigma}}{d\sigma}\right|_{\sigma=0} = \bE^{x}[\tau], 
\]
for $\eps=0.5$.

In order to obtain a reference solution with high accuracy we first compute $\psi_{\sigma}$ by discretizing the elliptic boundary value problem (\ref{FK}) based on a standard finite element discretization on a fine grid.
This is possible because the state space dimension in this guiding example is small but will not be possible in realistically high dimensions. The resulting reference solution for $\bE^{x}[\tau]$
is shown in the left panel of Figure \ref{fig:fmptUnc} below, along with the associated free energy $F_{\sigma}(x) = -\eps\log \psi_{\sigma}(x)$ in the right panel.

\begin{figure}
 \centering
 \includegraphics[width = 0.45\textwidth]{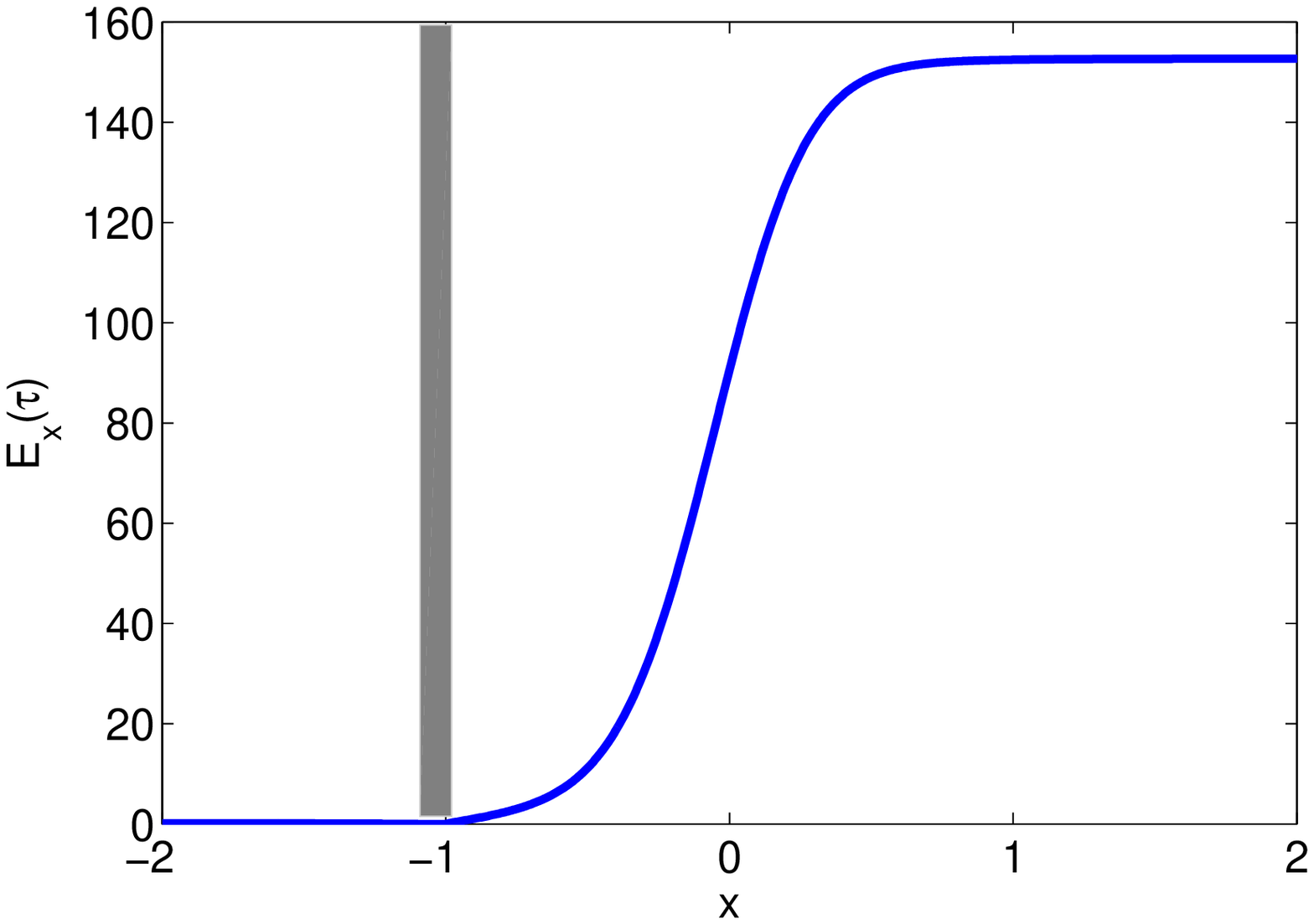}      
 \includegraphics[width = 0.45\textwidth]{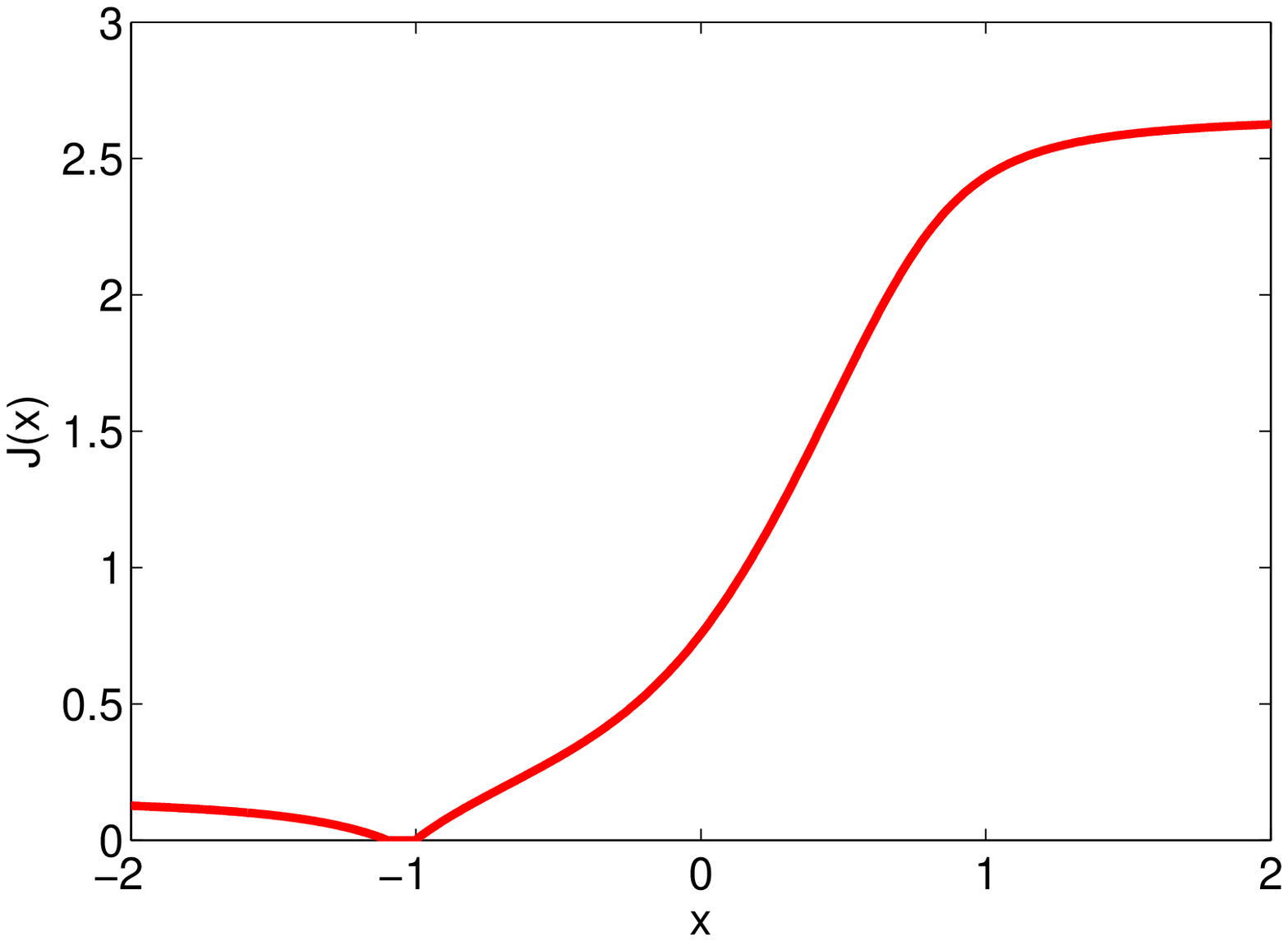}   
 \caption{Reference solution for the uncontrolled mean first passage time (left panel) and the related free energy $F_{\sigma}$ for $\sigma=1$ (right panel). Results based on  finite element discretization of (\ref{FK}) with high precision.}
 \label{fig:fmptUnc}
\end{figure}

\begin{figure}
 \centering
 \includegraphics[width = 0.45\textwidth]{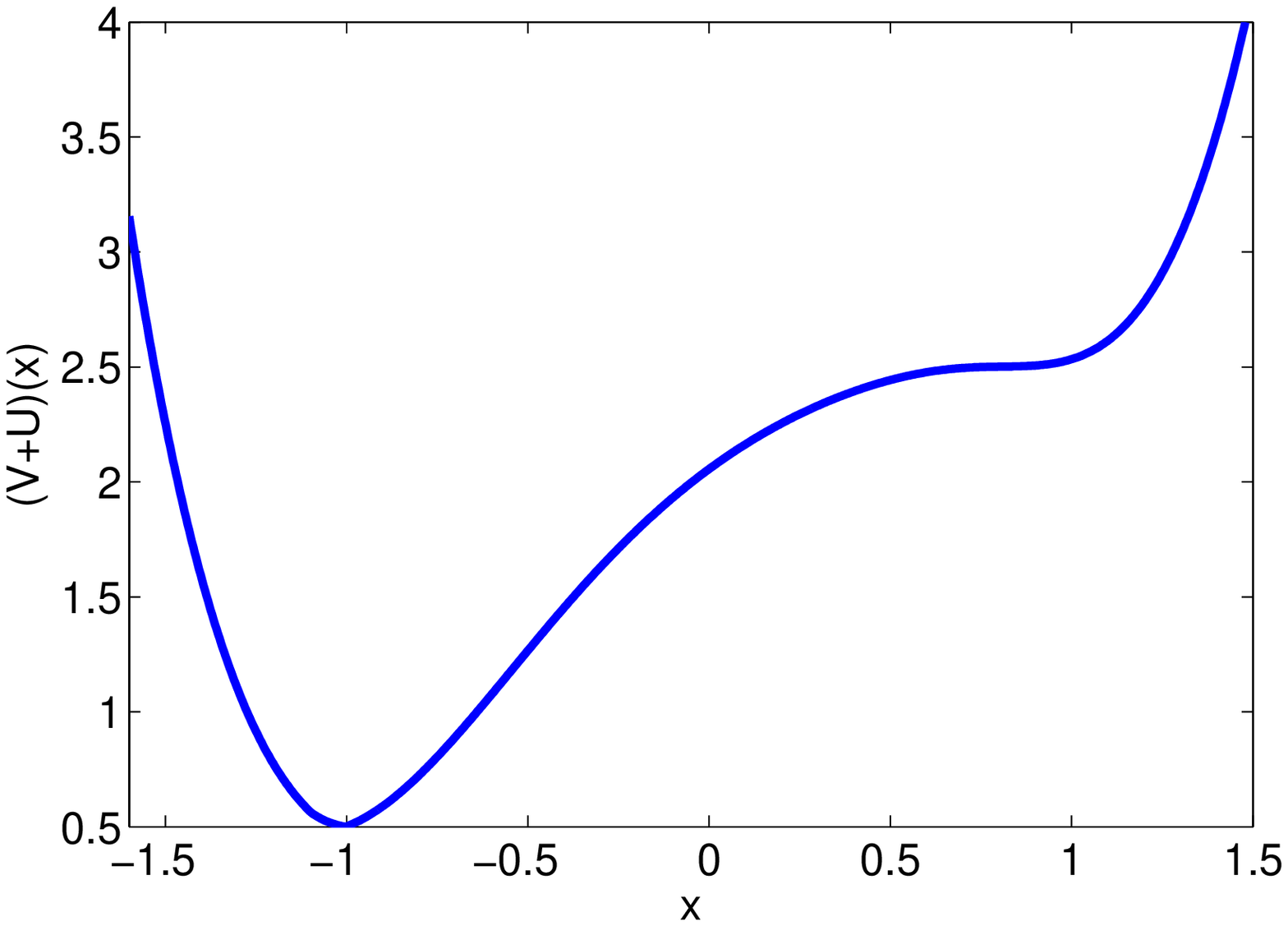}   
  \includegraphics[width = 0.45\textwidth]{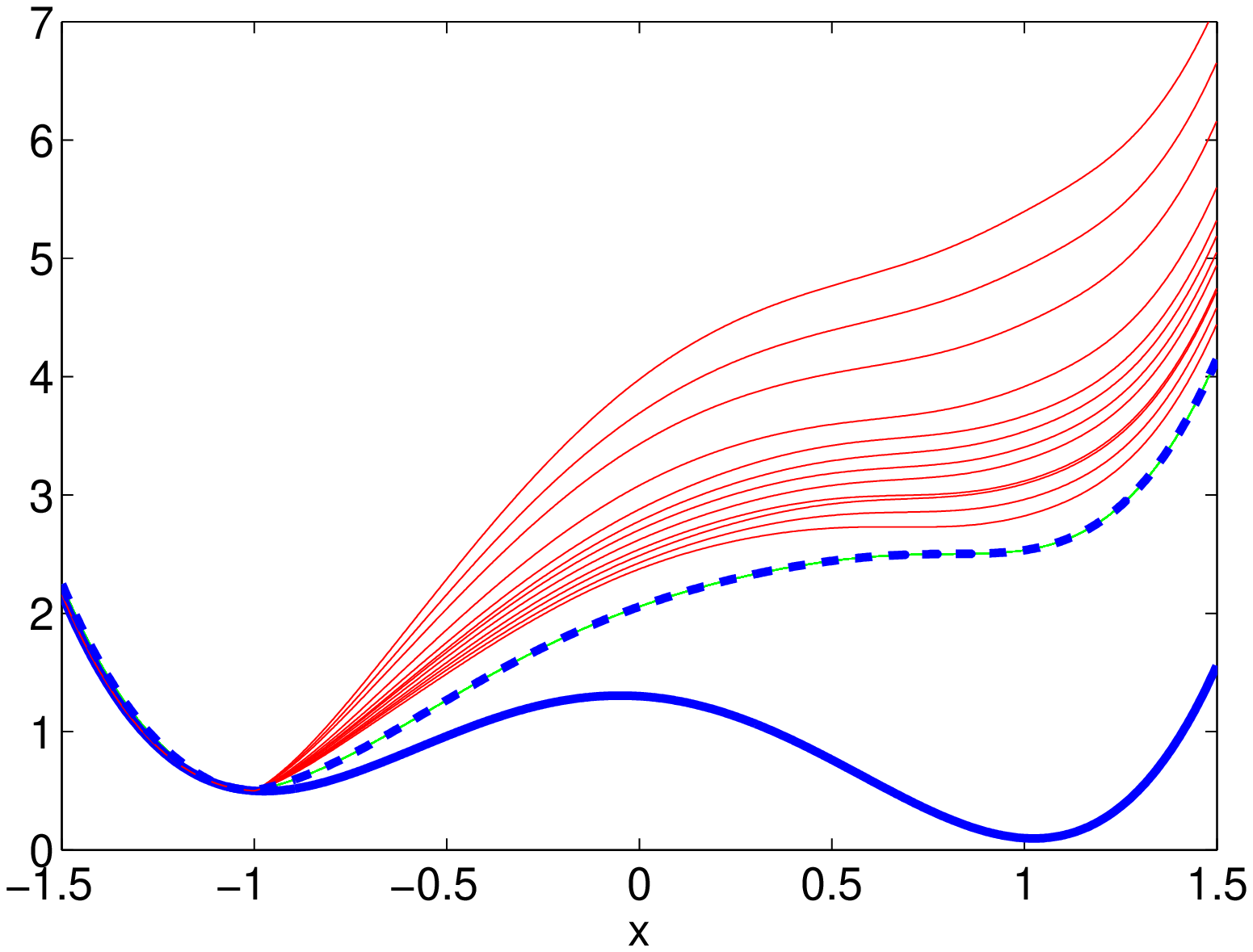}      
 \caption{Optimally tilted potential potential (left panel) and the first 11 iterates of the gradient descent (right panel). }
 \label{f:U}
\end{figure}

\begin{figure}
 \centering
 \includegraphics[width = 0.65\textwidth]{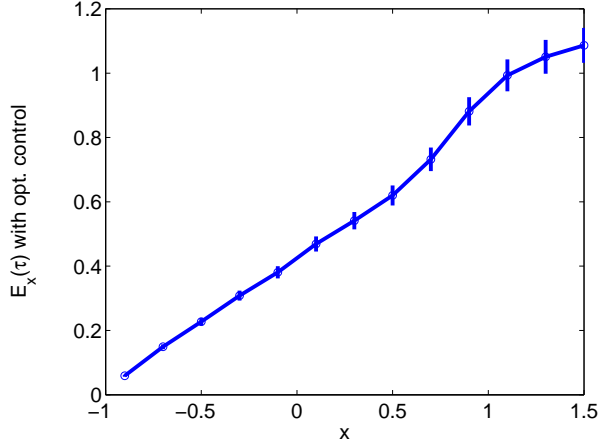}     
 \caption{Unbiased estimate of the first mean passage time, based on 2000 realizations of (\ref{tilted}) after Euler-discretization. The error bars indicate 95\% confidence intervals that were computed from the estimator's standard deviation.}
 \label{f:fmptC}
\end{figure}

An approximation of the free energy was then computed by the greedy milestoning / gradient descent algorithm described above that minimizes the cost functional (\ref{semicost}) in the coefficients $a=(a_{1},\ldots,a_{m})$. As scalar ansatz functions $v_j$ we chose $m=10$ Gaussians with width $0.1$ whose centers where uniformly spaced in the complement of $S$.  Once the minimization had been converged, the value function (free energy) and the resulting optimal control law were given by    
 \[
\tilde{F}= \sum_{j=1}^m a_j v_j(x)\,,\quad c^{*}(x)= \sum_{j=1}^m a_j b_j(x)\,,
\]
with $b_{j}=-\sqrt{2}\nabla v_{j}$. 
The result agree with the reference solution shown in Fig.~\ref{fig:fmptUnc} (deviations are of the order of the accurarcy threshold used in the gradient descent algorithm).
Figure \ref{f:U} shows the resulting optimally tilted potential $G=V+2\tilde{F}$, together with first few iteration steps of the gradient search. The mean first passage time of the tilted system
\begin{equation}\label{tilted} 
dX_{t} = - \nabla G(X_{s})dt + \sqrt{2\eps}dB_{s}\,,
\end{equation}
i.e., with $V$  in (\ref{sde}) replaced by the new potential $G$, is shown in Figure \ref{f:fmptC}.

As has been outlined above the algorithm only requires the computation of rather short trajectories since for all iterative potentials the mean first passage time is orders of magnitude smaller than for the original dynamics; the mean first passage time of the optimally tilted potential, e.g., is around 100 times smaller than originally.

\section{Conclusions and outlook}\label{sec:conclusion}
We have developed a simulation scheme for rare events that is based on an optimal change of measure that boils down to a logarithmic transformation of the path functional under consideration. The measure transformation turns the original exponential path functional into the functional of an optimal control problem that is linear in the observable and quadratic in the control variables. Although analytic solutions to the optimal control problem are available only in simple situations and computing the optimal change of measure may require to solve a possibly high-dimensional optimal control problem numerically, there is a considerable speed-up coming from (a) the fact that the functional is linear-quadratic and allows for the design of robust unbiased Monte-Carlo estimators and (b) the fact that events that were rare originally are no longer rare under the new probability measure. The gain in the numerical complexity requires that the optimal control problem can be solved efficiently, and, with the equivalence between path sampling and optimal control in hand, we have sketched a numerical algorithm for computing the optimal control that is based on an easy-to-implement inexact gradient descent that can be solved rather efficiently using milestoning. The algorithm was tested, computing the optimal feedback for the controlled passage between metastable sets in a double-well potential. Even though the numerical example that we presented is tiny on the scale of typical molecular dynamics applications, we emphasize that the minimization algorithm is independent of the dimension of the system and hence admits an easy generalization to more complicated systems; we refer to the rich literature on machine learning and queuing networks where various strategies for treating high dimensional systems have been developed (e.g., see \cite{deboer2005}). Finally we note that all ideas presented in this article can be readily extended to more complicated dynamics (e.g., degenerate diffusions with dissipation) and time-dependent path functionals (e.g., to simulate single-molecule experiments); it is even possible to consider situations where the exponential path functional involves additional control variables, in which case a logarithmic transformation leads to a game rather than an optimal control problem (cf.~\cite{latorre2010}). Further open issues are the deterministic limit of the stochastic control problem, the convergence analysis of the gradient descent and the rigorous analysis of fluctuations in systems under feedback control (cf.~\cite{ueda2010}).

\appendix

\section{Computational aspects}

In order to compute the gradient of (\ref{semicost}) \wrt to the unknown coefficients $a=(a_{1},\ldots,a_{m})$, it is convenient to discretize the path measure $Q=Q(a)$. To this end, let $0=t_{0}<t_{1}<\ldots<t_{N}=\tau$ be a set of time nodes with $h = t_{k+1}-t_{k}$ where we assume for the moment that $\tau<\infty$ is deterministic. Euler's method applied to 
\begin{equation*}
dX_{t} = \left(\sqrt{2} c(X_{t}) - \nabla V(X_{t})\right)dt + \sqrt{2\eps}\, dB_{t}\,.
\end{equation*}
gives 
\[
\tilde{X}_{k+1} = \tilde{X}_k + h\left(\sqrt{2}c(X_{k}) - \nabla V(\tilde{X}_{k})\right) + \sqrt{2 h \eps}\,\eta_{k+1}
\]
where the $\eta_k$ are i.i.d.~random variables that are normally distributed with mean zero and unit covariance. Since the $\eta_{k}$ are Gaussian, the density of the distribution $Q_{h}(a)$ of discrete paths $(\tilde{X}_{0},\ldots,\tilde{X}_{N})\subset\R^{n}$ conditional on $\tilde{X}_{0}=x_{0}$ is readily shown to be 
\begin{equation}\label{pathprob}
\varphi_{h}(x_{0},\ldots,x_{N};a) = (Z_{h}(a))^{-1} \exp\left(-S_{h}(x_{0},\ldots,x_{N};a)\right)
\end{equation}
with the discrete action
\begin{equation}\label{action}
S_{h} = \frac{h}{4\eps}\sum_{k=0}^{N-1}\left|\frac{x_{k+1}-x_{k}}{h}  + \nabla V(x_{k}) - \sqrt{2}c(x_{k}) \right|^{2} 
\end{equation}
and the normalization constant
\begin{equation}\label{partfct}
Z_{h} = \int_{\R^{n}\times\ldots\times\R^{n}} \exp\left(-S_{h}(x_{0},\ldots,x_{N};a)\right)dx_{1}\ldots dx_{N}\,.
\end{equation}
Computing the gradient of the discretized functional 
 \[
\tilde{I}_{h}(a;\tilde{X}_{0}) = \bE^{\tilde{X}_{0}}_{Q_{h}}\left[\sum_{k=0}^{N-1}\tilde{g}_{h}(\tilde{X}_{k}, c(\tilde{X}_{k}))\right]
\]
with 
 \[
\tilde{g}_{h}(x,c(x)) = h\left( f(x) + \frac{1}{2}\Big|\sum_{j}a_{j}b_{j}(x)\Big|^{2}\right)
\]
is now straightforward. Assuming that $\tilde{X}_{0}$ is independent of the control, we have 
\begin{equation*}
\frac{\partial \tilde{I}_{h}}{\partial a_j} =  \bE^{\tilde{X}_{0}}_{Q_{h}}\left[\sum_{k=0}^{N-1}\frac{\partial \tilde{g}_{h}}{\partial a_{j}} - \tilde{g}_{h}\left(\frac{\partial S_{h}}{\partial a_{j}}  + \frac{1}{Z_{h}}\frac{\partial Z_{h}}{\partial a_{j}}\right)  \right],
\end{equation*}
where both $\tilde{g}_{h}$ and $\partial \tilde{g}_{h}/\partial a_{j}$ are evaluated at $(x_{k},c(x_{k}))$. 
Specifically, 
\begin{eqnarray*}
\frac{\partial \tilde{g}_{h}}{\partial a_{j}} & = & h\, c(x_{k},t_{k}) b_j(x_{k})\\
\frac{\partial S_{h}}{\partial a_{j}} & = & -\frac{h}{\eps\sqrt{2}}\sum_{k=0}^{N-1}\left(\frac{x_{k+1}-x_{k}}{h}  + \nabla V(x_{k}) - \sqrt{2}c(x_{k})\right) b_{j}(x_{k})\\
\frac{\partial Z_{h}}{\partial a_{j}}& = & -Z_{h}\bE^{\tilde{X}_{0}}_{Q_{h}}\left[\frac{\partial S_{h}}{\partial a_{j}} \right].
\end{eqnarray*}
Together with the projection property of the conditional expectation this gives
\begin{equation}\label{gradI}
\frac{\partial \tilde{I}_{h}}{\partial a_j}  = \bE^{\tilde{X}_{0}}_{Q_{h}}\left[ \frac{\partial \tilde{g}_{h}}{\partial a_{j}}  \right] + \mathbf{C}^{\tilde{X}_{0}}_{Q_{h}}\left[\tilde{g}_{h},\,\frac{\partial S_{h}}{\partial a_{j}} \right] ,
\end{equation}
where $\mathbf{C}_{Q_{h}}$ denotes the covariance operator
\[
\mathbf{C}_{Q_{h}}[u,v]=\bE_{Q_{h}}[uv]-\bE_{Q_{h}}[u]\bE_{Q_{h}}[v].
\]

\paragraph{Inexact gradient}

We are interested in the situation when $\tau$ in (\ref{indefinite1}) is a random stopping time rather than a fixed time; otherwise the optimal control policy would be a function of time, i.e., $u_{t}=c(X_{t},t)$. But in case that $\tau$ is a first entry time of a set $S\subset\R^{n}$, this stopping time $\tau=\tau(c)$ will be a function of the control. Hence the derivative of the cost functional \wrt the unknown control coefficients $a_{j}$ would involve additional derivatives of $\tau$ or its time-discrete counterpart $N_{\tau}$; for example, for the discretized running cost this would result in an expression like  
\[
\frac{\partial}{\partial a_{j}} \sum_{k=0}^{N_{\tau}-1} \tilde{g}_{h}(x_{k},c(x_k)) = \tilde{g}_{h}(x_{N_{\tau}-1},c(x_{N_{\tau}-1}))\frac{\partial N_{\tau}}{\partial a_{j}} + \sum_{k=0}^{N_{\tau}-1} \frac{\partial \tilde{g}_{h}}{\partial a_{j}}  \]

In principle the dependence of the stopping time on the control variable can be made explicit in terms of the solution to an elliptic boundary value problem for $\tau$, yet it is unclear how terms such as $\partial N_{\tau}/\partial a_{j}$ can be handled numerically efficiently. 

In many cases the gradient descent will also converge even though the gradient $\nabla\tilde{I}$ is not exact, and it turns out that the boundary cost in the last equations is typically small compared to the accumulated cost. Ignoring the contribution from the boundary terms in the derivatives hence gives a gradient descent method with inaccurate gradient. In our numerical example where $f=\sigma$ is constant, the inexact gradient reads 
\begin{eqnarray*}
\frac{\partial \tilde{I}_{h}}{\partial a_j} & = &  h \bE^{\tilde{X}_{0}}_{Q_{h}}\left[\sum_{k=0}^{N_\tau-1}c(x_k) b_j(x_k)\right]\\
& & -\frac{h^{3/2}}{\eps \sqrt{2}} \mathbf{C}^{\tilde{X}_{0}}_{Q_{h}}\left[\sum_{k=0}^{N_\tau-1}(\sigma+\frac{1}{2}|c(X_k)|^{2}),\,\sum_{k=0}^{N_\tau-1} \eta_{k+1} b_{j}(x_k)\right]\\
& = &  h \bE^{\tilde{X}_{0}}_{Q_{h}}\left[\sum_{k=0}^{N_\tau-1}c(x_k) b_j(x_k)\right]\\
& & -\frac{h^{3/2}}{\eps \sqrt{2}}  \bE^{\tilde{X}_{0}}_{Q_{h}}\left[\left(\sum_{k=0}^{N_\tau-1}(\sigma+\frac{1}{2}|c(X_k)|^{2})\right)\sum_{k=0}^{N_\tau-1} \eta_{k+1} b_{j}(x_k)\right]\,,
\end{eqnarray*}
where $N_\tau=\lceil\tau/h\rceil$ is the discrete analog of the first hitting time (here $\lceil x\rceil$ is the nearest integer larger than $x$), and we used the fact that $b_{j}(X_{k})$ and $\eta_{k+1}$ are independent.

\end{document}